\documentclass[a4paper,12pt]{article}
\usepackage{CJK,CJKnumb}
\usepackage{indentfirst}
\usepackage{amssymb}
\usepackage{amsmath}
\usepackage[dvips]{graphicx}
\usepackage{amsthm}
\usepackage{amscd}
\usepackage{latexsym}
\usepackage{mathrsfs}

    \date{}
    \allowdisplaybreaks[4] \footskip=15pt

    \numberwithin{equation}{section} \theoremstyle{plain}
    \newtheorem*{thm*}{Main Theorem}

    \newtheorem*{cor*}{Corollary}
    
    \newtheorem*{lem*}{Lemma}
    
    \newtheorem*{prop*}{Proposition}
    
    \newtheorem*{rem*}{Remark}
    
    \newtheorem*{exa*}{Example}
    
    \newtheorem*{df*}{Definition}
    
    \newtheorem*{conj*}{Conjecture}
    \newtheorem*{ack*}{ACKNOWLEDGEMENTS}

    \newcommand{\pf}{\noindent\begin {proof}}
    \newcommand{\epf}{\end{proof}}

\begin{document}
\begin{center}
{\large\bf The structure of decomposable lattices determined \\ by
their prime ideals}\vskip 8mm

{\bf Xinmin Lu$^{*}$\,\,\,\,\,\,\, Dongsheng Liu} \vskip 1mm

{\small\it School of Science, Nanjing University of Science and
Technology,\\ Nanjing 210094, P.R. China}\vskip 2mm

{\bf Zhinan Qi} \vskip 1mm

{\small\it Department of Mathematics, Nanchang University,\\
Nanchang 330047, P.R. China}\vskip 2mm

{\bf Hourong Qin} \vskip 1mm

{\small\it Department of Mathematics, Nanjing University,\\
Nanjing 210093, P.R. China}\vskip 4mm

\begin{minipage}{130mm}
{\noindent\small\bf Abstract: }{\small  A distributive lattice $L$
with minimum element $0$ is called decomposable if $a$ and $b$ are
not comparable elements in $L$ then there exist
$\overline{a},\overline{b}\in L$ such that
$a=\overline{a}\vee(a\wedge b), b=\overline{b}\vee(a\wedge b)$ and
$\overline{a}\wedge \overline{b}=0$. The main purpose of this paper is to
 study the structure of decomposable lattices determined by
their prime ideals. The properties for five special decomposable lattices are derived.}  \vskip 2mm

\vskip 2mm

\vskip 2mm
\end{minipage}
\end{center}
$\\$

\footnotetext {$^{*}$ Corresponding author: School of Science,
Nanjing University of Science and Technology, Nanjing 210094, P.R.
China. E-mail: xmlu\_nanjing@hotmail.com}

\noindent{\bf 1. Introduction and main results}\vskip 4mm

Following [10], a decomposable lattice is a distributive lattice
$L$ with minimum element $0$ such that for any $a,b\in L$, if $a$
and $b$ are not comparable elements in $L$, then there exist
$\overline{a},\overline{b}\in L$ such that
$a=\overline{a}\vee(a\wedge b), b=\overline{b}\vee(a\wedge b)$ and
$\overline{a}\wedge \overline{b}=0$.  The idea of decomposable lattice is originated
from that of normal lattices and relatively normal lattices (see
e.g. [5,6,9,11,12,13]). We have described prime ideals, minimal prime ideals and
special ideals of a decomposable lattice in [10].

There are lots of decomposable lattices. In fact, it includes all
distributive and strongly projectable lattices, all the positive
cones of complete and compactly generated lattice-ordered groups,
and the lattices of ideals of some arithmetical rings. In [10] the
authors first established respectively a series of
characterizations of prime ideals, minimal prime ideals and
special ideals of a decomposable lattice and then investigated the
relationship among them. All these characterizations will be our
main technical tool for the further study of the structure of such
lattices. In the present paper, we shall apply the results in [10]
to study the structure of decomposable lattices determined by
their prime ideals. All results in this paper are purely
lattice-theoretic extension of some results of lattice-ordered
groups (see e.g. [2,3,4,7]).

Here is a brief outline   of the article. We simultaneously state
the main results.

In Section 2, we simply
review some basic definitions and introduce some notations for the
classes of decomposable lattices satisfy some special conditions.

In Section 3, we investigate decomposable lattices in which every
prime ideal contains at most $n$ minimal prime ideals. By using
the results in [10] and the pigeonhole principle, we shall prove
that the every prime ideal of a decomposable lattice $L$ contains
at most $n$ minimal prime ideals if and only if for any $n+1$
mutually disjoint elements $a_{1},a_{2},\cdots, a_{n},a_{n+1}$,
$L=a^{\perp}_{1}\vee a^{\perp}_{2}\vee\cdots\vee a^{\perp}_{n}\vee
a^{\perp}_{n+1}$.

In Section 4, we investigate decomposable lattices with basis and
prove that the following conditions are equivalent for a
decomposable lattice $L$: (1) $L$ has a basis; (2) for any $0<x\in
L$, there exists a basic element $a$ such that $x\geq a$; (3)
$P(L)$ is atomic; (4) for any $A\in P(L)\setminus \{L\}$,
$A=\bigcap\{P\in P(L)|\,\, P\supseteq A$, and $P$ is a maximal
polar ideal of $L \}$; (5) $\bigcap\{P\in P(L)|\,\, P$ is a
maximal polar ideal of $L\}=0$. As an application of this result,
we further investigate decomposable lattices with finite basis and
prove that the following conditions are equivalent for a
decomposable lattice $L$: (1) $L$ has a finite basis; (2) $P(L)$
is finite; (3) $P(L)$ satisfies $DCC$.

In Section 5, we investigate decomposable lattices with compact
property in the sense of Bigard-Conrad-Wolfenstein [1] and prove
that the following conditions are equivalent for a decomposable
lattice $L$: (1) $L$ is compact; (2) $L$ is discrete and each
minimal prime ideal of $L$ is a polar; (3) for any $M\in
MinSpe(L)$, there exists an atom $a$ of $L$ such that $a\not\in
M$; (4) each ultrafilter of $L$ is principal. This result is
purely lattice-theoretic extension of the corresponding result of
lattice-ordered groups. We apply the result to further investigate
the relationship between compact property and countably compact
property.

In Section 6, we investigate decomposable lattices $L$ in which
$V(L)$, $P(L)$ and $Ide(L)$ satisfy $DCC$, respectively, and prove
that $Ide(L)$ satisfies $DCC$ if and only if $V(L)$ and $P(L)$
satisfy $DCC$, respectively.

In last section, we investigate decomposable lattices in which
each nonzero element has only finitely many values and
decomposable lattices in which each disjoint subset with upper
bound is finite. Moreover, we also investigate consistency of
decomposable lattices  and establish a
simply connection between consistency and projectivity in the
category of decomposable lattices in which each nonzero element
has only finitely many values. \vskip 8mm

\noindent{\bf 2. Preliminaries and notations}\vskip 4mm

In this section, we simply review some basic definitions and some
well-known results. The readers are refereed to [8] for the general
theory of lattices.

Throughout this paper, we consider lattices $L$ with minimum
element $0$, denote by $\mathbb{DL}$ the class of decomposable
lattices and use "$\subset$" and "$\supset$" to denote proper
set-inclusion.

A lattice $L$ is called distributive if $a\wedge(b\vee c)=(a\wedge
b)\vee(a\wedge c)$ for any $a,b,c\in L$. A nonempty subset $I$ in
a lattice $L$ is called an ideal of $L$ if $a\vee b\in I$ for any
$a,b\in I$ and $a\geq x\in L$ implies that $x\in I$. We denote by
$Ide(L)$ the set of all ideals of $L$. In particular, if $a\in L$
then $(a]=\{x\in L|\,\, x\leq a\}$ is called the principal ideal
of $L$ generated by $a$. A direct computation shows that if $L\in
\mathbb{DL}$ then $Ide(L)$ is a distributive lattice by the rule:
$I\wedge J=I\cap J$ and $I\vee J=\{a\vee b|\,\, a\in I, b\in J\}$
for any $I,J\in Ide(L)$.

An ideal $P$ in a lattice $L$ is called prime if $P\neq L$ and
$a\wedge b\in P$ implies that either $a\in P$ or $b\in P$, where
$a,b\in L$. By Zorn's Lemma, each prime ideal contains a minimal
prime ideal. We denote by $Spe(L)$ and $MinSpe(L)$ respectively
the set of all prime ideals of $L$ and the set of all minimal
prime ideals of $L$.

Let $L$ be a lattice. For any $0<x\in L$, by Zorn's Lemma, there
exists a maximal ideal of $L$ with respect to not containing $x$,
denoted $M$, $M$ is called a regular ideal and is the value of
$x$. In general, $a$ need not have a unique value. We denote by
$Val(x)$ the set of all values of $x$. If $M$ is the unique value
of $x$, $M$ or $x$ is called special. We denote by $V(L)$ and
$S(L)$ respectively the set of all values of $L$ and the set of
all special values of $L$. Clearly, $S(L)\subseteq V(L)$. Observe
that the following conditions are equivalent: (1) $M\in V(L)$; (2)
$M$ is meet-irreducible, i.e., if
$\bigcap\limits_{\lambda\in\Lambda}I_{\lambda}=M$, where
$\{I_{\lambda}\}_{\lambda\in\Lambda}\subseteq Ide(L)$, then
$I_{\lambda}=M$ for some $\lambda$; (3) $M\subset
M^{*}=\bigcap\{I\in Ide(L)|\,\, I\supset M\}$; (4) $M\in Val(x)$,
where $x\in M^{*}\setminus M$.

For a lattice $L$ and $\emptyset\neq A\subseteq L$, we write
$A^{\perp}=\{x\in L|\,\, x\wedge a=0$ for any $a\in L \}$.
$A^{\perp}$ is called the polar of $A$, and define
$(A^{\perp})^{\perp}=A^{\perp\perp}$. $P\in Ide(L)$ is called
polar if $P=A^{\perp}$ for some $\emptyset\neq A\subseteq L$.
Clearly, $P\in Ide(L)$ is polar if and only if $P=P^{\perp\perp}$.
We denote by $P(L)$ the set of all polar ideals of $L$.

An element $a$ in a lattice $L$ is called a basic element if $a>0$
and $(a]$ is totally ordered. A nonempty subset
$\{a_{\lambda}\}_{\lambda\in \Lambda}$ of $L$ is called a basis if
this set is a maximal disjoint subset in $L$ and each element is a
basic element.

Let $L$ be a lattice and $\emptyset\neq F\subseteq L$. A nonempty
subset $F$ of $L$ is called a filter of $L$ if the following
conditions are satisfied: (1) $0\not\in F$; (2) for any $a,b\in
F$, $a\wedge b\in F$; (3) if $x\in L$ and $x\geq a\in F$ implies
$x\in F$. By Zorn's Lemma, each filter $F$ of $L$ must be
contained in a maximal filter $U$ of $L$, and $U$ is called an
ultrafilter of $L$. A filter $F$ of $L$ is called principal if
$F=\{x\in L|\,\, x\geq a\}$ for some $a\in L$.\vskip 2mm

In this article, $L$ will be always  a lattice unless otherwise
stated. For convenience, we  use the following notations to
denote classes of special lattices.\vskip 1mm

$\mathbb{A}=\{L|$\, every prime ideal of $L$ is minimal
$\}$.\vskip 1mm

$\mathbb{B}=\{L|$\, every prime ideal of $L$ contains a unique
minimal prime ideals $\}$.\vskip 1mm

$\mathbb{B}_{n}=\{L|$\, every prime ideal of $L$ contains at most
$n$ minimal prime ideals $\}$.\vskip 1mm

$\mathbb{B}_{\omega}=\{L|$\, every prime ideal of $L$ contains at
most finitely many minimal prime ideals$\}$.\vskip 1mm

$\mathbb{C}=\{L|$\, $L$ is compact $\}$.\vskip 1mm

$\mathbb{C}_{\omega}=\{L|$\, $L$ is countably compact $\}$.\vskip
1mm

$\mathbb{D}=\{L|$\, $Ide(L)$ satisfies $DCC \}$.\vskip 1mm

$\mathbb{E}=\{L|$\, $V(L)$ satisfies $DCC \}$.\vskip 1mm

$\mathbb{F}=\{L|$\, every disjoint subset of $L$ with upper bound
is finite $\}$.\vskip 1mm

$\mathbb{F}_{v}=\{L|$\, every nonzero element of $L$ has only
finitely many values $\}$.\vskip 1mm

$\mathbb{S}=\{L|$\, $L$ has a basis $\}$.\vskip 1mm

$\mathbb{S}_{\omega}=\{L|$\, $L$ has a finite basis $\}$. \vskip
1mm

$\mathbb{T}=\{L|$\, $L$ is projectable, i.e., for any $a\in L$,
$L=a^{\perp\perp}\vee a^{\perp}$ $\}$.\vskip 8mm

\noindent{\bf 3. $\mathbb{B}$ and $\mathbb{B}_{n}$}\vskip 4mm

In this section, we shall investigate decomposable lattices in
which every prime ideal contains at most $n$ minimal prime ideals.
By using the results in [10] and the pigeonhole principle, we
shall establish explicit characterizations for the class of such
lattices.

First, we need the following two lemmas ([10], Lemma 4.2 and Lemma
5.7). \vskip 2mm

\noindent{\bf Lemma 3.1. }Let $L\in \mathbb{DL}$. If $P\in Spe(L)$
then\vskip 2mm
\begin{center}
$\bigcup\{a^{\perp}|\,\, a\in L\setminus P\}=\bigcap\{M\in
MinSpe(L)|\,\, M\subseteq P\}$.\vskip 3mm
\end{center}

\noindent{\bf Lemma 3.2. }Let $L\in \mathbb{DL}$. If
$Q_{1},Q_{2},\cdots, Q_{n}$ are mutually incomparable prime ideals
of $L$ and $a\not\in Q_{i}$ for $i=1,2,\cdots,n$, then there exist
$a_{i}\in (\bigcap\limits_{j\neq i}Q_{j}) \setminus Q_{i}$ such
that $0<a_{i}<a$ for $i=1,2,\cdots,n$ and $a_{i}\wedge a_{j}=0$
for $i\neq j$.\vskip 2mm

We now state and prove the main result of this section.\vskip 2mm

\noindent{\bf Theorem 3.3. }Let $L\in \mathbb{DL}$. The following
conditions are equivalent:

(1) $L\in \mathbb{B}_{n}$.

(2) For any $n+1$ distinct minimal prime ideals
$M_{1},M_{2},\cdots, M_{n},M_{n+1}$ of $L$,\vskip 2mm
\begin{center}
$L=M_{1}\vee M_{2}\vee\cdots\vee M_{n}\vee M_{n+1}$.\vskip 2mm
\end{center}

(3) For any $n+1$ mutually incomparable values
$Q_{1},Q_{2},\cdots, Q_{n},Q_{n+1}$ of $L$,\vskip 2mm
\begin{center}
$L=Q_{1}\vee Q_{2}\vee\cdots\vee Q_{n}\vee Q_{n+1}$.\vskip 2mm
\end{center}

(4) For any $n+1$ mutually incomparable prime ideals
$P_{1},P_{2},\cdots, P_{n},P_{n+1}$ of $L$,\vskip 2mm
\begin{center}
$L=P_{1}\vee P_{2}\vee\cdots\vee P_{n}\vee P_{n+1}$.\vskip 2mm
\end{center}

(5) For any $n+1$ mutually disjoint elements $a_{1},a_{2},\cdots,
a_{n},a_{n+1}$,\vskip 2mm
\begin{center}
$L=a^{\perp}_{1}\vee a^{\perp}_{2}\vee\cdots\vee a^{\perp}_{n}\vee
a^{\perp}_{n+1}$.\vskip 3mm
\end{center}

\noindent{\bf Proof.
}(1)$\Rightarrow$(2)$\Rightarrow$(3)$\Rightarrow$(4)$\Rightarrow$(1)
is clear. It suffices to show (1)$\Leftrightarrow$(5).\vskip 2mm

(1)$\Rightarrow$(5) Assume that there exist $n+1$ mutually
disjoint elements $a_{1},a_{2},\cdots, a_{n},a_{n+1}$ in $L$ such
that $\bigvee\limits_{i=1}^{n+1}a^{\perp}_{i}\subset L$. Pick
$x\in L\setminus (\bigvee\limits_{i=1}^{n+1}a^{\perp}_{i})$. Then
there exists some $M\in Val(x)$ such that $M\supseteq
\bigvee\limits_{i=1}^{n+1}a^{\perp}_{i}$. By Lemma 3.1, we
have\vskip 2mm
\begin{center}
$\bigcap\{P\in MinSpe(L)|\,\, P\subseteq
M\}=\bigcup\{h^{\perp}|\,\, h\not\in M\}$.\vskip 2mm
\end{center}

Now, write $V_{0}=\bigcup\{h^{\perp}|\,\, h\not\in M\}$. We claim
that $a_{i}\not\in V_{0}$ for any $i$ ($i=1,2,\cdots,n,n+1$).
Assume that $a_{i}\in V_{0}$ for some $i$. Then $a_{i}\in
h^{\perp}$ for some $h\not\in M$. Thus $h\in
a^{\perp}_{i}\subseteq M$, a contradiction. Since $L\in
\mathbb{B}_{n}$, $M$ contains at most $n$ minimal prime ideals,
write $Q_{1},Q_{2},\cdots, Q_{n}$. Then\vskip 2mm
\begin{center}
$V_{0}=\bigcap\limits_{j=1}^{n}Q_{j}=\bigcup\{h^{\perp}|\,\,
h\not\in M\}$.\vskip 2mm
\end{center}
\noindent By the pigeonhole principle, there exists a minimal
prime ideal, denoted $Q_{i}$, contained in $M$, which does not
contain two of the elements of the set $\{a_{1},a_{2},\cdots,
a_{n},a_{n+1}\}$. Since $Q_{i}$ is prime, this is not impossible.
So $L=a^{\perp}_{1}\vee a^{\perp}_{2}\vee\cdots\vee
a^{\perp}_{n}\vee a^{\perp}_{n+1}$.

(5)$\Rightarrow$(1) Assume that there exists some $Q\in Spe(L)$
such that $Q$ contains $n+1$ distinct minimal prime ideals of $L$,
write $Q_{1},Q_{2},\cdots, Q_{n},Q_{n+1}$. Clearly, they are
mutually incomparable. So, by Lemma 3.2, there exist\vskip 2mm
\begin{center}
$a_{i}\in (\bigcap\limits_{1\leq j\neq i\leq n+1}Q_{j}) \setminus
Q_{i}$, where $i,j=1,2,\cdots,n,n+1$,\vskip 2mm
\end{center}
\noindent such that $a_{i}\wedge a_{j}=0$ for $i\neq j$. Now, for
any $i$ ($i=1,2,\cdots,n,n+1$), since $a_{i}\not\in Q_{i}$,
$a^{\perp}_{i}\subseteq Q_{i}$. So\vskip 2mm
\begin{center}
$Q\supseteq Q_{1}\vee Q_{2}\vee\cdots\vee Q_{n}\vee
Q_{n+1}\supseteq a^{\perp}_{1}\vee a^{\perp}_{2}\vee\cdots\vee
a^{\perp}_{n}\vee a^{\perp}_{n+1}=L$,\vskip 2mm
\end{center}
\noindent which is a contradiction. Therefore $L\in
\mathbb{B}_{n}$.$\hfill\Box$\vskip 2mm

As a direct result of Theorem 3.3, we have\vskip 2mm

\noindent{\bf Corollary 3.4. }Let $L\in \mathbb{DL}$. If for any
$M\in MinSpe(L)$, $L=M\vee M^{\perp}$, then $L\in
\mathbb{B}$.\vskip 8mm

\noindent{\bf 4. $\mathbb{S}$ and $\mathbb{S}_{\omega}$}\vskip 4mm

In this section, we investigate decomposable lattices with basis
or finite basis. We shall establish a series of characterizations
for them.

Recall that an element $a$ in a lattice $L$ is called a basic
element if $a>0$ and $(a]$ is totally ordered. A nonempty subset
$\{a_{\lambda}\}_{\lambda\in \Lambda}$ of $L$ is called a basis if
this set is a maximal disjoint subset in $L$ and each element is a
basic element.

The following two lemmas are well known ([10], Theorem 4.6 and
Theorem 4.3). \vskip 2mm

\noindent{\bf Lemma 4.1. }Let $L\in \mathbb{DL}$ and $0\neq I\in
Ide(L)$. The following conditions are equivalent:

(1) $I$ is totally ordered.

(2) For any $0<a\in I$, $a^{\perp}=I^{\perp}$.

(3) $I^{\perp}\in Spe(L)$.

(4) $I^{\perp}\in MinSpe(L)$.

(5) $I^{\perp\perp}$ is a maximal totally ordered ideal of $L$.

(6) $I^{\perp\perp}$ is a minimal polar ideal of $L$.

(7) $I^{\perp}$ is a maximal polar ideal of $L$.

(8) For any $0<a\in I$, $a$ is special.\vskip 2mm

\noindent{\bf Lemma 4.2. }Let $L\in \mathbb{DL}$ and $P\in
Spe(L)$. The following conditions are equivalent:

(1) $P\in MinSpe(L)$.

(2) $P=\bigcup\{a^{\perp}|\,\, a\not\in P\}$.

(3) For any $x\in P$, $x^{\perp}\not\subseteq P$.\vskip 2mm

We shall first apply Lemma 4.1 and Lemma 4.2 to establish
characterizations of decomposable lattices with basis.\vskip 2mm

\noindent{\bf Theorem 4.3. }Let $L\in \mathbb{DL}$. The following
conditions are equivalent:

(1) $L\in \mathbb{S}$.

(2) For any $0<x\in L$, there exists a basic element $a$ such that
$x\geq a$.

(3) $P(L)$ is atomic, i.e., for any $A\in P(L)\setminus \{0,L\}$,
there exists a minimal polar ideal $B$ of $L$ such that
$A\supseteq B$.

(4) For any $A\in P(L)\setminus \{L\}$, $A=\bigcap\{P\in P(L)|\,\,
P\supseteq A$, and $P$ is a maximal polar ideal of $L \}$.

(5) $\bigcap\{P\in P(L)|\,\, P$ is a maximal polar ideal of
$L\}=0$.\vskip 2mm

\noindent{\bf Proof. }(1)$\Rightarrow$(2) Let
$\{a_{\lambda}\}_{\lambda\in\Lambda}$ be a basis of $L$. Now, for
any $0<x\in L$, since $\{a_{\lambda}\}_{\lambda\in\Lambda}$ is a
maximal disjoint subset in $L$, there exists some
$\lambda\in\Lambda$ such that $x\wedge a_{\lambda}>0$. Since
$a_{\lambda}$ is a basic element, $x\wedge a_{\lambda}$ is clearly
a basic element and $x\geq x\wedge a_{\lambda}$.

(2)$\Rightarrow$(3) Given any $A\in P(L)\setminus \{0,L\}$, pick
$0<x\in A$. By (2), there exists a basic element $a$ such that
$x\geq a$. Then $a^{\perp\perp}\subseteq x^{\perp\perp}\subseteq
A$. By Lemma 4.1, $a^{\perp\perp}$ is a minimal polar ideal of
$L$. So $P(L)$ is atomic.

(3)$\Rightarrow$(4) Since the map $P\rightarrow P^{\perp}$ for any
$P\in P(L)$ is a dual isomorphism of lattices, by (3), for any
$A\in P(L)\setminus \{L\}$, there exists a maximal polar ideal $P$
of $L$ such that $P\supseteq A$. Consider the set\vskip 2mm
\begin{center}
$\Omega=\{P\in P(L)|\,\, P\supseteq A$, and $P$ is a maximal polar
ideal of $L \}$.\vskip 2mm
\end{center}
\noindent Clearly, $A\subseteq\bigcap\limits_{P\in \Omega}P$. If
$0<x\not\in A=A^{\perp\perp}$, then there exists some $0<a\in
A^{\perp}$ such that $x\wedge a>0$. By (3), there exists a maximal
polar ideal $P$ such that $(x\wedge a)^{\perp}\subseteq P$. So
$A=A^{\perp\perp}\subseteq a^{\perp}\subseteq (x\wedge
a)^{\perp}\subseteq P$. By Lemma 4.1, $P$ is a maximal polar ideal
implies that $P\in MinSpe(L)$. So, by Lemma 4.2, $x\wedge a\not\in
P$, and hence $x\not\in P$. Therefore $A=\bigcap\limits_{P\in
\Omega}P$.\vskip 1mm

(4)$\Rightarrow$(5) Suppose that $\bigcap\{P\in P(L)|\,\, P$ is a
maximal polar ideal of $L\}\neq 0$. Pick $0<a\in \bigcap\{P\in
P(L)|\,\, P$ is a maximal polar ideal of $L\}$. By (4),
$a^{\perp}\subseteq P$ for some maximal polar ideal $P$ of $L$.
Again, by Lemma 4.1, $P\in MinSpe(L)$. But $a\in P$ and
$a^{\perp}\subseteq P$, which contradicts Lemma 4.2. So
$\bigcap\{P\in P(L)|\,\, P$ is a maximal polar ideal of $L\}=0$.

(5)$\Rightarrow$(1) Let
$\bigcap\limits_{\lambda\in\Lambda}P_{\lambda}=0$, where each
$P_{\lambda}$ is a maximal polar ideal of $L$. Then each
$P^{\perp}_{\lambda}$ is a minimal polar ideal of $L$. By Lemma
4.1, $P^{\perp}_{\lambda}$ is totally ordered. Now, pick
$0<a_{\lambda}\in P^{\perp}_{\lambda}$ for any $\lambda\in
\Lambda$. Clearly, each $a_{\lambda}$ is a basic element of $L$.
Set $A=\{a_{\lambda}|\,\, \lambda\in\Lambda\}$. We shall show that
$A$ is a basis of $L$. For any $\alpha,\beta\in \Lambda$ with
$\alpha\neq\beta$, $P^{\perp}_{\alpha}$ and $P^{\perp}_{\beta}$
are both minimal polar ideals of $L$, $a_{\alpha}\wedge
a_{\beta}\in P^{\perp}_{\alpha}\cap P^{\perp}_{\beta}=0$. In
addition, if $x\wedge a_{\lambda}=0$ for any $\lambda\in \Lambda$,
where $x\in L$, then $x\in
a^{\perp}_{\lambda}=P^{\perp\perp}_{\lambda}=P_{\lambda}$ for any
$\lambda\in \Lambda$. Thus $x\in
\bigcap\limits_{\lambda\in\Lambda}P_{\lambda}=0$. So
$A=\{a_{\lambda}|\,\, \lambda\in\Lambda\}$ is a basis of $L$.
Therefore $L\in \mathbb{S}$.$\hfill\Box$\vskip 2mm

As an application of Theorem 4.3, we have\vskip 2mm

\noindent{\bf Corollary 4.4. }Let $L\in \mathbb{DL}$. If every
minimal prime ideal of $L$ is a polar ideal, i.e.,
$MinSpe(L)\subseteq P(L)$, then $L\in \mathbb{S}$.\vskip 2mm

\noindent{\bf Proof. }Given any $M\in MinSpe(L)$, there exists
$\emptyset\neq A\subseteq L$ such that $M=A^{\perp}$. Clearly,
$A\neq \{0\}$. Pick $0<a\in A$. Then $M=A^{\perp}\subseteq
a^{\perp}$.

First, we claim that $M= a^{\perp}$. Assume that $a^{\perp}\supset
M$. Then $a\not\in M$. Pick $0<b\in a^{\perp}\setminus M$. Since
$a\wedge b=0\in M$ and $M$ is prime, this means that either $a\in
M$ or $b\in M$, a contradiction. So $M=a^{\perp}$.

Second, we show that $a$ is a basic element. Otherwise, there
exist $0<a_{1},a_{2}<a$ such that $a_{1}\wedge a_{2}=0\in
M=a^{\perp}$. So either $a_{1}\in a^{\perp}$ or $a_{2}\in
a^{\perp}$, which implies that either $a_{1}=0$ or $a_{2}=0$, a
contradiction.

Finally, we show that $L$ has a basis. Assume that $L$ has no
basis. By Theorem 4.3, there exists $0<x\in L$ such that $x$ does
not exceed any basic elements. Let $Q_{x}$ be a value of $x$.
Since every prime ideal of $L$ contains at least a minimal prime
ideal. Without loss of generality, suppose that $M\subseteq
Q_{x}$. Clearly, $x\not\in M=a^{\perp}$, so $x\wedge a>0$. Notice
that $x\wedge a$ is a basic element and $x\geq x\wedge a$, a
contradiction. Therefore $L\in \mathbb{S}$.$\hfill\Box$ \vskip 2mm

For a decomposable lattice $L$ and $M\in V(L)$, $M$ is called
essential if there exists $0<x\in L$ such that for any
$G_{\lambda}\in Val(x)$, $G_{\lambda}\subseteq M$. We denote by
$E(L)$ the set of all essential values of $L$. Clearly,
$S(L)\subseteq E(L)$. Write\vskip 2mm
\begin{center}
$Rad(L)=\bigcap E(L)$,\vskip 2mm
\end{center}
\noindent and $Rad(L)$ is called the radical of $L$.

In the following, we shall use Theorem 4.3 to establish a
connection between decomposable lattices with basis and
$Rad(L)=0$.\vskip 2mm

\noindent{\bf Corollary 4.5. }Let $L\in \mathbb{DL}$ and $L\in
\mathbb{B}_{\omega}$. The following conditions are equivalent:

(1) $L\in \mathbb{S}$.

(2) $\bigcap S(L)=0$.

(3) $Rad(L)=0$.\vskip 2mm

\noindent{\bf Proof. }(1)$\Rightarrow$(2) Let $\{s_{\lambda}|\,\,
\lambda\in\Lambda\}$ be a basis of $L$. Assume that $\bigcap
S(L)\neq 0$. Pick $0<x\in \bigcap S(L)$. By Theorem 4.3, there
exists some $\lambda\in\Lambda$ such that $x\geq s_{\lambda}$.
Since $s_{\lambda}$ is special, let $Q_{\lambda}$ be the unique
value of $s_{\lambda}$, then $s_{\lambda}\not\in Q_{\lambda}$, so
that $s_{\lambda}\not\in \bigcap S(L)$, which contradicts the fact
that $x\in \bigcap S(L)$ implies $s_{\lambda}\in \bigcap S(L)$.

(2)$\Rightarrow$(3) Since $S(L)\subseteq E(L)$, $Rad(L)\subseteq
\bigcap S(L)=0$. So $Rad(L)=0$.

(3)$\Rightarrow$(1) Given any $0<a\in L$, since $Rad(L)=0$, there
exists some $M\in E(L)$ such that $a\not\in M$. So there exists
$Q\in Val(a)$ such that $Q\supseteq M$. Since $M\in E(L)$, $Q\in
E(L)$.

Now, we claim that there exists a basic element $s$ of $L$ such
that $a\geq s$. Since $Q$ is essential, there exists $0<f\in L$
such that all the values of $f$ are contained in $Q$. Set
$g=a\wedge f$. Clearly, all the values of $g$ are also contained
in $Q$. For convenience, we may suppose that $f<a$. Since $L\in
\mathbb{B}_{\omega}$, there exists a positive integer $k$ such
that $Q$ contains at most $k$ minimal prime ideals. Now, if $f$
does not exceed a basic element, then there exists a disjoint
subset $\{f_{1},f_{2},\cdots,f_{m}\}$ of $L$ with upper bound $f$
and satisfies $m>k$. Let $Q_{i}$ be a value of $f_{i}$ for
$i=1,2,\cdots,m$. Clearly, $Q_{i}\parallel Q_{j}$ for $i\neq j$
and each $Q_{i}\subseteq Q$, a contradiction. So $a$ must exceed a
basic element of $L$. Therefore $L\in
\mathbb{S}$.$\hfill\Box$\vskip 2mm

In order to establish characterizations of decomposable lattices
with finite basis, we need the following lemma.\vskip 2mm

\noindent{\bf Lemma 4.6. }Let $L\in \mathbb{DL}$ and let
$\{a_{1},a_{2},\cdots,a_{n}\}$ be a basis of $L$ and
$a^{\perp\perp}_{i}=A_{i}$ for $1\leq i\leq n$. Then
$(\bigvee\limits_{i\in \bigtriangleup}A_{i})^{\perp}=
(\bigvee\limits_{i\in
N\setminus\bigtriangleup}A_{i})^{\perp\perp}$, where
$N=\{1,2,\cdots,n\}$ and $\bigtriangleup\subseteq N$.\vskip 3mm

\noindent{\bf Proof. }We divide the proof into two steps.\vskip
2mm

{\bf Step 1}. If $\bigtriangleup=\emptyset$ then
$\bigvee\limits_{i\in \bigtriangleup}A_{i}=0$, we are done.\vskip
3mm

{\bf Step 2}. If $\bigtriangleup\neq\emptyset$ then $a_{i}\wedge
a_{j}=0$ for any $i\neq j$. Thus $a_{i}\in a^{\perp}_{j}$, and
hence $a^{\perp\perp}_{j}\subseteq a^{\perp}_{i}$, so that
$a^{\perp\perp}_{i}\cap a^{\perp\perp}_{j}=0$, i.e., $A_{i}\cap
A_{j}=0$. So\vskip 3mm
\begin{center}
$\bigvee\limits_{i\in \bigtriangleup}A_{i}\subseteq
(\bigvee\limits_{i\in
N\setminus\bigtriangleup}A_{i})^{\perp}\Rightarrow(\bigvee\limits_{i\in
N\setminus\bigtriangleup}A_{i})^{\perp\perp}\subseteq(\bigvee\limits_{i\in
\bigtriangleup}A_{i})^{\perp}$.\vskip 3mm
\end{center}

On the other hand, $\{a_{1},a_{2},\cdots,a_{n}\}$ is a basis of
$L$, so $(\bigvee\limits_{i\in N}A_{i})^{\perp}=0$. Thus\vskip 3mm
\begin{center}
$(\bigvee\limits_{i\in
\bigtriangleup}A_{i})^{\perp}\cap(\bigvee\limits_{i\in
N\setminus\bigtriangleup}A_{i})^{\perp}=0$,\vskip 3mm
\end{center}
\noindent so that $(\bigvee\limits_{i\in
\bigtriangleup}A_{i})^{\perp}\subseteq (\bigvee\limits_{i\in
N\setminus\bigtriangleup}A_{i})^{\perp\perp}$. Therefore
$(\bigvee\limits_{i\in \bigtriangleup}A_{i})^{\perp}=
(\bigvee\limits_{i\in
N\setminus\bigtriangleup}A_{i})^{\perp\perp}$.$\hfill\Box$ \vskip
4mm

\noindent{\bf Theorem 4.7. }Let $L\in \mathbb{DL}$. The following
conditions are equivalent:

(1) $L\in \mathbb{S}_{\omega}$.

(2) $P(L)$ is finite.

(3) $P(L)$ satisfies $DCC$.

\noindent{\bf Proof. }(1)$\Rightarrow$(2) Let
$\{a_{1},a_{2},\cdots,a_{n}\}$ be a finite basis of $L$. Set
$A_{i}=a^{\perp\perp}_{i}$ for $1\leq i\leq n$. Then each $A_{i}$
is a minimal polar ideal of $L$ by Lemma 4.1. So, for any $P\in
P(L)$, either $P\cap A_{i}=0$ or $P\cap A_{i}=A_{i}$. Set
$\bigtriangleup=\{1,2,\cdots,n\}$, and\vskip 2mm
\begin{center}
$\bigtriangleup_{1}=\{i\in \bigtriangleup|\,\, P\cap
A_{i}=0\}$\,\, and \,\, $\bigtriangleup_{2}=\{i\in
\bigtriangleup|\,\, P\cap A_{i}=A_{i}\}$.\vskip 2mm
\end{center}
\noindent So\vskip 2mm
\begin{center}
$\bigvee\limits_{i\in \bigtriangleup_{1}}A_{i}\subseteq
P^{\perp}$\,\,\, and \,\,\, $\bigvee\limits_{i\in
\bigtriangleup_{2}}A_{i}\subseteq P$.\vskip 2mm
\end{center}
\noindent Then\vskip 2mm
\begin{center}
$\bigvee\limits_{i\in \bigtriangleup_{1}}A_{i}\subseteq
P^{\perp}\subseteq (\bigvee\limits_{i\in
\bigtriangleup_{2}}A_{i})^{\perp}\Rightarrow(\bigvee\limits_{i\in
\bigtriangleup_{2}}A_{i})^{\perp\perp}\subseteq
P^{\perp\perp}\subseteq (\bigvee\limits_{i\in
\bigtriangleup_{2}}A_{i})^{\perp}$.\vskip 3mm
\end{center}
\noindent By Lemma 4.6, $P=P^{\perp\perp}=(\bigvee\limits_{i\in
\bigtriangleup_{2}}A_{i})^{\perp\perp}$. Therefore $P(L)$ is
finite.\vskip 1mm

(2)$\Rightarrow$(3) is clear.

(3)$\Rightarrow$(1) We first show that $L\in \mathbb{S}$.
Otherwise, there exists $0<a\in L$ such that $a$ does not exceed
any basic element. Hence there exist $0<a_{0},b_{0}<a$ such that
$a_{0}\wedge b_{0}=0$. For $a_{0}$, $a_{0}$ does not exceed any
basic element. Hence there exist $0<a_{1},b_{1}<a_{0}$ such that
$a_{1}\wedge b_{1}=0$. Continuing this process, we can obtain an
infinite descending chain of $P(L)$ as follows:\vskip 2mm
\begin{center}
$a^{\perp}_{0}\supset a^{\perp}_{1}\supset\cdots\supset
a^{\perp}_{n}\supset\cdots$,\vskip 2mm
\end{center}
which contradicts the fact that $P(L)$ satisfies $DCC$. So $L$
must has a basis. Now, let $\{s_{\lambda}|\,\,
\lambda\in\Lambda\}$ be a basis of $L$. Assume that
$|\Lambda|=\infty$. Then we can similarly obtain an infinite
descending chain of $P(L)$ as follows:\vskip 2mm
\begin{center}
$s^{\perp}_{1}\supset (s_{1}\vee
s_{2})^{\perp}\supset\cdots\supset (s_{1}\vee s_{2}\vee\cdots\vee
s_{n})^{\perp}\supset\cdots$,\vskip 2mm
\end{center}
\noindent a contradiction. Therefore $L\in
\mathbb{S}_{\omega}$.$\hfill\Box$\vskip 2mm

Recall that an element $u$ in a lattice $L$ is called a unit if
$u>0$ and $u\wedge x>0$ for any $0<x\in L$. Recall also that $I\in
Ide(L)$ is called large if $I\cap J\neq 0$ for any $0\neq J\in
Ide(L)$. As an application of Theorem 4.7, we have\vskip 2mm

\noindent{\bf Corollary 4.8. }Let $L\in \mathbb{DL}$ and $L\in
\mathbb{B}_{\omega}$. The following conditions are equivalent:

(1) $L\in \mathbb{S}_{\omega}$.

(2) There exists a large ideal $I$ of $L$ with the form
$I=\bigvee\limits_{i=1}^{n}(a_{i}]$, where each $(a_{i}]$ is
totally ordered.

(3) $L$ has a unit $u$ and $u$ has only finitely many
values.\vskip 2mm

\noindent{\bf Proof. }(1)$\Rightarrow$(2) Let
$\{a_{1},a_{2},\cdots,a_{n}\}$ be a finite basis of $L$. A direct
computation shows that $I=\bigvee\limits_{i=1}^{n}(a_{i}]$ is a
large ideal of $L$, and each $(a_{i}]$ is clearly totally ordered.

(2)$\Rightarrow$(3) For any $0<x\in L$, since $I$ is large and
$Ide(L)$ is a distributive lattice, we then have \vskip 2mm
\begin{center}
$0\neq (x]\cap I=(x]\cap(\bigvee\limits_{i=1}^{n}(
a_{i}])=\bigvee\limits_{i=1}^{n}((x]\cap(a_{i}])$.\vskip 2mm
\end{center}
\noindent So there exists some $i$ ($1\leq i\leq n$) such that $(
x]\cap( a_{i}]\neq 0$, i.e., $0\neq x\wedge a_{i}\leq x$. By
Theorem 4.3, $\{a_{1},a_{2},\cdots,a_{n}\}$ is a finite basis of
$L$. Now, let $u=a_{1}\vee a_{2}\vee\cdots\vee a_{n}$. Clearly,
$u$ is a unit of $L$. Since each $a_{i}$ is a basic element, each
$a_{i}$ must be special. Let $Q_{i}$ be the unique value of
$a_{i}$. Since $\{a_{1},a_{2},\cdots,a_{n}\}$ are mutually
disjoint, we have\vskip 2mm
\begin{center}
$Val(u)=\bigcup\limits_{i=1}^{n}Val(a_{i})=\{Q_{1},Q_{2},\cdots,Q_{n}\}$.\vskip
2mm
\end{center}
\noindent So $u$ has only finitely many values.

(3)$\Rightarrow$(1) Let $\{Q_{1},Q_{2},\cdots,Q_{k}\}$ be the set
of all values of $u$.

We first show that $L\in \mathbb{F}$. Otherwise, there exists
$0<f\in L$ and an infinite disjoint subset of $L$ with upper bound
$f$, write $\{a_{i}\in L|\,\, i\in I, |I|=\infty \}$. Notice that
$u$ is a unit. Set $b_{i}=u\wedge a_{i}$ for any $i\in I$. Then
$\{b_{i}\in L|\,\, i\in I, |I|=\infty \}$ is an infinite disjoint
subset of $L$ with upper bound $u$. Now, let $M_{i}$ be a value of
$b_{i}$ for any $i\in I$. Since $u\not\in M_{i}$ for any $i\in I$,
$M_{i}\subseteq Q_{j}$ for some $j$ ($j=1,2,\cdots,k$). Notice
that $I$ is infinite, so that there exists some $Q_{j}$ contains
an infinite number of $M_{i}$, which contradicts $L\in
\mathbb{B}_{\omega}$.

Since $L\in \mathbb{B}_{\omega}$, let $Q_{i}$ contains $n_{i}$
minimal prime ideals, and set $m=\max\{n_{i}|\,\, 1\leq i\leq
k\}$. We shall show that if the number of basic elements of $L$ is
$n$, then $n\leq mk+1$. Otherwise, there exists a disjoint subset
$\{x_{1},x_{2},\cdots,x_{t}\}$ of $L$ such that $t>mk+1$. Set
$y_{i}=x_{i}\wedge u$ for $i=1,2,\cdots,t$. Then
$\{y_{1},y_{2},\cdots,y_{t}\}$ is also a mutually disjoint subset
of $L$ with upper bound $u$. Repeating the above process, we shall
obtain that there exists some $Q_{i}$ contains at least $m+1$
minimal prime ideals, a contradiction. Therefore $L\in
\mathbb{S}_{\omega}$.$\hfill\Box$\vskip 8mm

\noindent{\bf 5. $\mathbb{C}$ and $\mathbb{C}_{\omega}$}\vskip 4mm

In this section, we shall first study the structure of
decomposable lattices with compact property and then investigate
the relationship between compact property and countably compact
property.

Recall that a lattice $L$ is called compact if
$\{a_{\lambda}\}_{\lambda\in\Lambda}$ is a nonempty subset of $L$
and $\bigwedge\limits_{\lambda\in\Lambda}a_{\lambda}=0$ then there
exists a finite subset $\{a_{i}\}^{n}_{i=1}$ of
$\{a_{\lambda}\}_{\lambda\in\Lambda}$ such that
$\bigwedge\limits_{i=1}^{n}a_{i}=0$. Recall also that a lattice
$L$ is called discrete if every nonzero element of $L$ exceeds an
atom.\vskip 2mm

\noindent{\bf Theorem 5.1. }Let $L\in \mathbb{DL}$. The following
conditions are equivalent:

(1) $L\in \mathbb{C}$.

(2) $L$ is discrete and each minimal prime ideal of $L$ is a
polar.

(3) For any $M\in MinSpe(L)$, there exists an atom $a$ of $L$ such
that $a\not\in M$.

(4) Each ultrafilter of $L$ is principal.\vskip 2mm

\noindent{\bf Proof. }(1)$\Rightarrow$(2) We first show that $L$
is discrete. Given any $0<x\in L$, let
$\{a_{\lambda}\}_{\lambda\in\Lambda}$ be a maximal chain of $L$
containing $x$. If
$\bigwedge\limits_{\lambda\in\Lambda}a_{\lambda}=0$, then since
$L$ is compact, there exist a finite subset $\{a_{i}\}^{n}_{i=1}$
of $\{a_{\lambda}\}_{\lambda\in\Lambda}$ such that
$\bigwedge\limits_{i=1}^{n}a_{i}=0$. Notice that
$a_{1},a_{2},\cdots,a_{n}$ are mutually comparable and each
$a_{\lambda}>0$, this is clearly impossible. So, if we set
$a=\bigwedge\limits_{\lambda\in\Lambda}a_{\lambda}$, then $a$ is
an atom in $L$ and $x\geq a$. Thus $L$ is discrete.

We next show that every minimal prime ideal of $L$ is a polar.
Given any $M\in MinSpe(L)$, $M=\bigcup\{a^{\perp}|\,\, a\in
L\setminus M\}$ by Lemma 4.2. Set $K=L\setminus M$. Notice that
$L\in \mathbb{DL}$ and $M$ is prime, which implies that $K$ is a
chain of $L$. So $K$ has an atom $a$ such that $M=a^{\perp}$.

(2)$\Rightarrow$(3) For any $M\in MinSpe(L)$, by (2), $M\in P(L)$.
So $M=A^{\perp}$ for some $\emptyset \neq A\subseteq L$. Now pick
$0<x\in A$. Since $L$ is discrete, there exists an atom $a$ in $L$
such that $x\geq a$. We claim that $a\not\in M$. Otherwise, $a\in
M=A^{\perp}$ implies $a\in A\cap A^{\perp}=0$, a contradiction.

(3)$\Rightarrow$(4) Let $K$ be an ultrafilter of $L$. Then
$M=\bigcup\{a^{\perp}|\,\, a\in K\}$ is a minimal prime ideal of
$L$. By (3), there exists an atom $a\not\in M$ such that
$a^{\perp}\subseteq M$. Notice that $a$ is an atom implies that
$a^{\perp}\in MinSpe(L)$ by Lemma 4.1, and hence $M=a^{\perp}$.

Now, it suffices to show that $K=\{x\in L|\,\, x\geq a\}$. For any
$x\in K$, then $x^{\perp}\subseteq M=a^{\perp}$. So $x\wedge a>0$
for any $x\in K$. Since $a$ is an atom, $x\wedge a=a$. Thus $x\geq
a$. Conversely, given any $x\in L$, if $x\geq a$ then since
$a\not\in M$, this means $x\not\in M$, so that $x\in K$. Thus
$K=\{x\in L|\,\, x\geq a\}$. So $K$ is principal.

(4)$\Rightarrow$(1) Let $\{a_{\lambda}\}_{\lambda\in\Lambda}$ be a
nonempty subset of $L$. Suppose that for any finite subset
$\{a_{i}\}^{n}_{i=1}$ of $\{a_{\lambda}\}_{\lambda\in\Lambda}$,
$\bigwedge\limits_{i=1}^{n}a_{i}\neq 0$. Set
$a=\bigwedge\limits_{i=1}^{n}a_{i}$. Let $Q\in Val(a)$ and $M\in
MinSpe(L)$ be such that $M\subseteq Q$. Then $a\not\in M$. Since
the set $\{a_{\lambda}\}_{\lambda\in\Lambda}$ must contained in an
ultrafilter of $L$, write $K$. By (4), $K=\{x\in L|\,\, x\geq b\}$
for some $0<b\in L$. So
$\bigwedge\limits_{\lambda\in\Lambda}a_{\lambda}\geq b>0$.
Therefore $L\in \mathbb{C}$.$\hfill\Box$\vskip 2mm

We now apply Theorem 5.1 to establish the relationship between
compact property and countably compact property.

Let us recall that a lattice $L$ is called countably compact if
$\{a_{i}\}_{i=1}^{\infty}$ is a subset of $L$ and
$\bigwedge\limits_{i=1}^{\infty}a_{i}=0$ then there exists a
positive integer $n$ such that
$\bigwedge\limits_{i=1}^{n}a_{i}=0$.\vskip 2mm

\noindent{\bf Theorem 5.2. }Let $L\in \mathbb{DL}$. If
$MinSpe(L)\subseteq P(L)$ then following conditions are
equivalent:

(1) $L\in \mathbb{C}_{\omega}$.

(2) Every totally ordered ideal of $L$ is countably compact.\vskip
2mm

\noindent{\bf Proof. }(1)$\Rightarrow$(2) is clear.

(2)$\Rightarrow$(1) By way of contradiction. Assume that there
exists a countable subset $\{a_{i}\}_{i=1}^{\infty}$ of $L$ such
that $\bigwedge\limits_{i=1}^{\infty}a_{i}=0$, but
$\bigwedge\limits_{i=1}^{n}a_{i}\neq 0$ for any positive integers
$n$.

First, by Corollary 4.4, $L$ has a basis. Let $\{a_{\lambda}|\,\,
\lambda\in\Lambda\}$ be a basis of $L$. We claim that for any
$\lambda\in\Lambda$ there exists some positive integer $i$ such
that $a_{\lambda}\wedge a_{i}=0$. Otherwise, there exists some
$\mu\in \Lambda$ such that $a_{\mu}\wedge a_{i}>0$ for any
positive integers $i$. Since\vskip 2mm
\begin{center}
$\bigwedge\limits_{i=1}^{\infty}(a_{\mu}\wedge
a_{i})=a_{\mu}\wedge(\bigwedge\limits_{i=1}^{\infty}a_{i})=0$,\vskip
2mm
\end{center}
\noindent and\vskip 2mm
\begin{center}
$\{a_{\mu}\wedge a_{i}|\,\, i=1,2,\cdots \}\subseteq (
a_{\mu}]$,\vskip 2mm
\end{center}
\noindent by (2), there exists a positive integer $n$ such that
$\bigwedge\limits_{i=1}^{n}(a_{\mu}\wedge a_{i})=0$. Notice that
each $a_{\mu}\wedge a_{i}>0$ and $\{a_{\mu}\wedge a_{i}|\,\,
i=1,2,\cdots,n \}$ is a finite chain, this is clearly impossible.

Second, since $\{a_{i}\}_{i=1}^{\infty}$ is a $\wedge$-semilattice
of $L$, this means that there exists an ultrafilter $U$ of $L$
such that $\{a_{i}\}_{i=1}^{\infty}\subseteq U$. So
$Q=\bigcup\{a^{\perp}|\,\, a\in U\}$ is a minimal prime ideal of
$L$.

Finally, for any $\lambda\in\Lambda$, using the above result,
there exists some positive integer $i$ such that
$a_{\lambda}\wedge a_{i}=0\in Q$ and $a_{i}\not\in Q$, then
$a_{\lambda}\in Q$ for any $\lambda\in\Lambda$. So
$\{a_{\lambda}|\,\, \lambda\in\Lambda\}\subseteq Q$. But, by
hypothesis, $Q=A^{\perp}$ for some $\emptyset\neq A\subseteq L$,
which implies that $Q=a^{\perp}_{\lambda}$ for some basic element
$a_{\lambda}$ in $L$. This is clearly impossible. Therefore $L\in
\mathbb{C}_{\omega}$.$\hfill\Box$\vskip 8mm

\noindent{\bf 6. $\mathbb{D}$ and $\mathbb{E}$}\vskip 4mm

In this section, we shall investigate decomposable lattices $L$ in
which $Ide(L)$ satisfies $DCC$ and prove that $\mathbb{D}=
\mathbb{E}\cap \mathbb{S}_{\omega}$.

By a direct computation, we have\vskip 2mm

\noindent{\bf Lemma 6.1. }Let $L\in \mathbb{DL}$. Then $L\in
\mathbb{E}$ if and only if $Spe(L)=V(L)$.\vskip 2mm

\noindent{\bf Theorem 6.2. }Let $L\in \mathbb{DL}$. The following
conditions are equivalent:

(1) $L\in \mathbb{D}$.

(2) $V(L)$ and $P(L)$ satisfy $DCC$ respectively.

(3) $L\in \mathbb{E}\cap \mathbb{S}_{\omega}$.\vskip 2mm

\noindent{\bf Proof. }(1)$\Rightarrow$(2) is automatic and
(2)$\Leftrightarrow$(3) is clear by Lemma 6.1 and Theorem 4.7. It
suffices to show (2)$\Rightarrow$(1).

By way of contradiction. Assume that there exists an infinite
descending chain of $Ide(L)$, as follows:\vskip 2mm
\begin{center}
$I_{1}\supset I_{2}\supset\cdots\supset I_{n}\supset\cdots$.\vskip
2mm
\end{center}
\noindent Pick $a_{i}\in I_{i}\setminus I_{i+1}$ for
$i=1,2,\cdots,n,\cdots$. Let $Q_{i}$ be a value of $a_{i}$ with
$Q_{i}\supseteq I_{i+1}$.

Now, let $\{a_{1},a_{2},\cdots,a_{n}\}$ be a finite basis of $L$.
We claim that $L$ has only $n$ minimal prime ideals
$a^{\perp}_{1},a^{\perp}_{2},\cdots,a^{\perp}_{n}$. In fact, it
suffices to show that for any $P\in MinSpe(L)$, there exists some
$i$ ($1\leq i\leq n$) such that $P=a^{\perp}_{i}$. Suppose that
$a^{\perp}_{i}\not\subseteq P$ for any $i$. Pick $b_{i}\in
a^{\perp}_{i}\setminus P$ for $i=1,2,\cdots,n$, and set
$b=b_{1}\wedge b_{2}\wedge\cdots\wedge b_{n}$. Then $0<b\not\in P$
and $b\wedge a_{i}=0$ for $i=1,2,\cdots,n$, which contradicts the
fact that $\{a_{1},a_{2},\cdots,a_{n}\}$ is a basis of $L$. So $L$
has only $n$ minimal prime ideals, written
$P_{1},P_{2},\cdots,P_{n}$. Notice that\vskip 2mm
\begin{center}
$I_{2}\subseteq Q_{1}, I_{3}\subseteq Q_{2}, \cdots,
I_{n+1}\subseteq Q_{n}, \cdots$.\vskip 2mm
\end{center}
\noindent Then $L$ has at least infinite many distinct values
$Q_{1},Q_{2},\cdots, Q_{n},\cdots$.

For $P_{1}$, since $L\in \mathbb{E}$, there exists a finite subset
of the set $\{Q_{1},Q_{2},\cdots, Q_{n},\cdots\}$ containing
$P_{1}$ which is a proper descending chain of $V(L)$. Now, we omit
this subset from the set $\{Q_{1},Q_{2},\cdots, Q_{n},\cdots\}$.

For $P_{2}$, similarly, there also exists a finite subset of the
set $\{Q_{1},Q_{2},\cdots, Q_{n},\cdots\}$ containing $P_{2}$
which is a proper descending chain of $V(L)$. We also omit this
subset from the set $\{Q_{1},Q_{2},\cdots, Q_{n},\cdots\}$.

Continuing this process, finally, for $P_{n}$, there also exists a
finite subset of the set $\{Q_{1},Q_{2},\cdots, Q_{n},\cdots\}$
containing $P_{n}$ which is a proper descending chain of $V(L)$.
We similarly omit this subset from the set $\{Q_{1},Q_{2},\cdots,
Q_{n},\cdots\}$.

Notice that the set $\{Q_{1},Q_{2},\cdots, Q_{n},\cdots\}$ is
infinite, the remains are also an infinite subset of
$\{Q_{1},Q_{2},\cdots, Q_{n},\cdots\}$, and each of which does not
contain any one of minimal prime ideals
$P_{1},P_{2},\cdots,P_{n}$. This is clearly impossible. So $L\in
\mathbb{D}$.$\hfill\Box$\vskip 8mm

\noindent{\bf 7. $\mathbb{F}_{v}$ and $\mathbb{F}$}\vskip 4mm

In this section, we investigate decomposable lattices in which
each nonzero element has only finitely many values and
decomposable lattices in which each disjoint subset with upper
bound is finite.

The following is well known ([10], Theorem 5.9).\vskip 2mm

\noindent{\bf Lemma 7.1. }Let $L\in \mathbb{DL}$. The following
conditions are equivalent:

(1) $L\in \mathbb{F}_{v}$.

(2) For any $0<a\in L$, $a=a_{1}\vee a_{2}\vee\cdots\vee a_{n}$,
where $a_{i}\wedge a_{j}=0$ for $i\neq j$ and each $a_{i}$ is
special.\vskip 2mm

\noindent{\bf Theorem 7.2. }Let $L\in \mathbb{DL}$. If $L\in
\mathbb{F}$, then $L\in \mathbb{A}$ and $V(L)=S(L)$. \vskip 2mm

\noindent{\bf Proof. }We first claim that $L$ has a basis.
Otherwise, there exists $0<a\in L$ and $0<a_{1}, a_{2}<a$ such
that $a_{1}\wedge a_{2}=0$. $a$ is not a basic element implies
that $a_{1}$ is not also a basic element, then there exist
$0<a_{11}, a_{12}<a_{1}$ such that $a_{11}\wedge a_{12}=0$.
Continuing this process, we can obtain an infinite disjoint subset
$\{x_{1},x_{2},\cdots,x_{n},\cdots\}$ of $L$ with upper bound $a$,
which contradicts $L\in \mathbb{F}$. So $L$ has a basis.

Now, let $\{a_{\lambda}|\,\, \lambda\in\Lambda\}$ be a basis of
$L$. We further claim that there exists $\lambda\in\Lambda$ such
that $a^{\perp}_{\lambda}\subseteq P$ for any $P\in MinSpe(L)$.
Suppose that $a^{\perp}_{\lambda}\not\subseteq P$ for any
$\lambda\in\Lambda$. We shall divide the proof into two
steps:\vskip 2mm

{\bf Step 1}. If $|\Lambda|=\infty$, then pick a fixed
$\lambda_{1}\in \Lambda$. Since
$a^{\perp}_{\lambda_{1}}\not\subseteq P$, we may further pick
$0<b_{1}\in a^{\perp}_{\lambda_{1}}\setminus P$. Since $L\in
\mathbb{F}$, there exists a finite subset of $\Lambda$, write
$\{\lambda_{2},\lambda_{3},\cdots,\lambda_{n}\}$, such that
$b_{1}\wedge a_{\lambda_{i}}>0$ for $i=2,3,\cdots,n$, and
$b_{1}\wedge a_{\lambda}=0$ for any $\lambda\in \Lambda\setminus
\{\lambda_{2},\lambda_{3},\cdots,\lambda_{n}\}$. Now, pick
$b_{i}\in a^{\perp}_{\lambda_{i}}\setminus P$ for
$i=2,3,\cdots,n$, and set $b=b_{1}\wedge b_{2}\wedge\cdots\wedge
b_{n}$. Then $0<b\not\in P$ and $b\wedge a_{\lambda}=0$ for any
$\lambda\in \Lambda$, which contradicts the fact that
$\{a_{\lambda}|\,\, \lambda\in\Lambda\}$ is a basis of $L$.\vskip
2mm

{\bf Step 2}. If $|\Lambda|=k<\infty$, then pick $b_{i}\in
a^{\perp}_{i}\setminus P$ for $i=1,2,\cdots,k$, and set
$b=b_{1}\wedge b_{2}\wedge\cdots\wedge b_{k}$. Then $0<b\not\in P$
and $b\wedge a_{i}=0$ for $i=1,2,\cdots,k$, which also contradicts
the fact that $\{a_{1},a_{2},\cdots,a_{k}\}$ is a basis of $L$.

In view of the above arguments, there exists some
$\lambda\in\Lambda$ such that $a^{\perp}_{\lambda}\subseteq P$. By
Lemma 4.1, $a^{\perp}_{\lambda}\in MinSpe(L)$, so that
$P=a^{\perp}_{\lambda}$.

We now show that $L\in \mathbb{A}$. Given any $M\in Spe(L)$, there
exists some $P\in MinSpe(L)$ such that $M\supseteq P=a^{\perp}$.
By Lemma 4.1, $a$ is a basic element, and hence $a^{\perp}$ is a
maximal poplar. So $M\supseteq a^{\perp}$, which implies that
$M^{\perp\perp}\supseteq a^{\perp}$, so that
$a^{\perp}=M^{\perp\perp}$. So $M=a^{\perp}\in MinSpe(L)$.
Therefore $L\in \mathbb{A}$.

Finally, we show that $S(L)=V(L)$. For any $Q\in V(L)$, using the
above result, $Q$ is a minimal prime ideal of $L$ and
$Q=s^{\perp}$, where $s$ is a basic element of $L$. So $Q$ is the
unique value of $s$, and hence $S(L)=V(L)$.$\hfill\Box$\vskip 2mm

\noindent{\bf Theorem 7.3. }Let $L\in \mathbb{DL}$. If $A\in
\mathbb{A}\cap\mathbb{F}_{v}$ then $L\in \mathbb{F}$. \vskip 2mm

\noindent{\bf Proof. }Assume that $L\not\in \mathbb{F}$. Then
there exists some $0<a\in L$ and an infinite disjoint subset
$\{a_{\lambda}|\,\, \lambda\in \Lambda\}$ of $L$ with upper bound
$a$. Now, let $Q_{\lambda}$ be a value of $a_{\lambda}$ for any
$\lambda\in \Lambda$. Since $L\in \mathbb{F}_{v}$, $a$ has only
finitely many values, write $Q_{1},Q_{2},\cdots,Q_{n}$. Since
$a_{\lambda}\not\in Q_{\lambda}$ for any $\lambda\in \Lambda$
implies that $a\not\in Q_{\lambda}$ for any $\lambda\in \Lambda$.
This means that for any $\lambda\in \Lambda$,
$Q_{\lambda}\subseteq Q_{i}$ for some $i$ ($i=1,2,\cdots, n$).
Since $L\in \mathbb{A}$, we get that for any $\lambda\in \Lambda$,
$Q_{\lambda}=Q_{i}$ for some $i$ ($i=1,2,\cdots, n$). Notice that
$\Lambda$ is infinite, this is impossible. So $L\in
\mathbb{F}$.$\hfill\Box$\vskip 2mm

At the end of this paper, we investigate consistency of
decomposable lattices in the sense of the following
definition.\vskip 2mm

\noindent{\bf Definition 7.4. }Let $L$ be a lattice. For any
$0<x\in L$, we denote by $v(x)$ the cardinal number of the set of
all values of $x$, i.e., $v(x)=|Val(x)|$. A lattice $L$ is called
consistent if  $x\leq y$ then $v(x)\leq v(y)$, where $0<x,y\in L$.
\vskip 2mm

According to Definition 7.4, one can obtain that for a
decomposable lattice $L$, if $L$ is consistent then $a$ is special
if and only if it is a basic element. In order to investigate the
structure of decomposable lattices with consistency, let us recall
that a lattice $L$ is called projectable if for any $a\in L$,
$L=a^{\perp\perp}\vee a^{\perp}$.\vskip 2mm

\noindent{\bf Theorem 7.5. }Let $L\in \mathbb{DL}$. If $L\in
\mathbb{F}_{v}$ then the following conditions are equivalent:

(1) $L\in \mathbb{T}$.

(2) $L\in \mathbb{B}$.

(2) $L$ is consistent.\vskip 2mm

\noindent{\bf Proof. }(1)$\Rightarrow$(2) Suppose that $L\not\in
\mathbb{B}$. Then there exist some $P\in Spe(L)$ and
$M_{1},M_{2}\in MinSpe(L)$ with $M_{1}\neq M_{2}$ such that
$M_{1}\vee M_{2}\subseteq P$. Pick $a\in M_{1}\setminus M_{2}$. By
Lemma 4.2, $a^{\perp\perp}\subseteq M_{1}$ and $a^{\perp}\subseteq
M_{2}$. Then\vskip 2mm
\begin{center}
$L=a^{\perp\perp}\vee a^{\perp}\subseteq M_{1}\vee M_{2}\subseteq
P$,\vskip 2mm
\end{center}
\noindent a contradiction. Thus $L\in \mathbb{B}$.

(2)$\Rightarrow$(3) Given $x,y\in L$ with $0<x\leq y$, let $Q_{x}$
be a value of $x$. Then $x\not\in Q_{x}$ implies that $y\not\in
Q_{x}$, and hence there exists some $Q_{y}\in Val(y)$ such that
$Q_{x}\subseteq Q_{y}$. In order to show that $v(x)\leq v(y)$, it
suffices to show that if $Q_{1}$ and $Q_{2}$ are two distinct
values of $x$, then $Q_{1}$ and $Q_{2}$ can not be contained in
the same value $Q$ of $y$. Otherwise, since $L\in \mathbb{B}$,
$L=Q_{1}\vee Q_{2}\subseteq Q$, a contradiction. Therefore $L$ is
consistent.

(3)$\Rightarrow$(1) Since $L$ is consistent, each special element
in $L$ is a basic element. Again, since $L\in \mathbb{F}_{v}$, by
Lemma 7.1, $L$ has a basis. Now, let $\{a_{\lambda}|\,\,
\lambda\in\Lambda\}$ be a basis of $L$. Then, by Lemma 4.1,
$\{a^{\perp\perp}_{\lambda}|\,\, \lambda\in\Lambda\}$ is a set of
maximal totally ordered ideals of $L$.\vskip 1mm

{\bf Claim 1}. If $D$ is a maximal totally ordered ideal of $L$,
then $D=a^{\perp\perp}_{\lambda_{0}}$ for some $\lambda_{0}\in
\Lambda$.\vskip 1mm

Notice that $\{a_{\lambda}|\,\, \lambda\in\Lambda\}$ is a basis of
$L$. We can obtain that $D\cap a^{\perp\perp}_{\lambda_{0}}\neq 0$
for some $\lambda_{0}\in \Lambda$. Otherwise, $D\cap
a^{\perp\perp}_{\lambda}=0$ for any $\lambda\in\Lambda$ implies
that $d\wedge a_{\lambda}=0$ for any $\lambda\in\Lambda$ and
$0<d\in D$, which contradicts the fact that $\{a_{\lambda}|\,\,
\lambda\in\Lambda\}$ is a basis of $L$. So
$D=a^{\perp\perp}_{\lambda_{0}}$ for some $\lambda_{0}\in
\Lambda$.\vskip 1mm

{\bf Claim 2}. For any $0<g\in L$, $g^{\perp\perp}=
\bigvee\limits_{\lambda_{i}\in\Lambda_{1}}a^{\perp\perp}_{\lambda_{i}}$,
where $\Lambda_{1}$ is a finite subset of $\Lambda$.\vskip 2mm

Since $L\in \mathbb{F}_{v}$, by Lemma 7.1, $g=g_{1}\vee
g_{2}\vee\cdots\vee g_{n}$, where $g_{i}\wedge g_{j}=0$ for $i\neq
j$ and each $g_{i}$ is special. Moreover, each $g_{i}$ is a basic
element. Using Claim 1, each $g_{i}\in
g^{\perp\perp}_{i}=a^{\perp\perp}_{\lambda_{i}}$ for some
$\lambda_{i}\in\Lambda$, where $i=1,2,\cdots,n$. Set\vskip 2mm
\begin{center}
$\Lambda_{1}=\{\lambda_{1},\lambda_{2},\cdots,\lambda_{n}\}$ and
$\Lambda_{2}=\Lambda\setminus \Lambda_{1}$.\vskip 2mm
\end{center}
\noindent Clearly, $g^{\perp\perp}=
\bigvee\limits_{\lambda_{i}\in\Lambda_{1}}a^{\perp\perp}_{\lambda_{i}}$.\vskip
3mm

{\bf Claim 3}. For any $0<g\in L$,
$g^{\perp}=\bigvee\limits_{\lambda\in\Lambda_{2}}a^{\perp\perp}_{\lambda}$.\vskip
3mm

For any $\lambda\in \Lambda_{2}$, $a_{\lambda}\wedge
a_{\lambda_{i}}=0$ for any $\lambda_{i}\in\Lambda_{1}$, which
implies that $a_{\lambda}\wedge g_{i}=0$ for $i=1,2,\cdots,n$. So
$a_{\lambda}\wedge g=0$.  Then $a_{\lambda}\in g^{\perp}$ for any
$\lambda\in \Lambda_{2}$, so that
$a^{\perp\perp}_{\lambda}\subseteq g^{\perp}$ for any $\lambda\in
\Lambda_{2}$. So
$\bigvee\limits_{\lambda\in\Lambda_{2}}a^{\perp\perp}_{\lambda}\subseteq
g^{\perp}$.\vskip 1mm

Conversely, for any $0<h\in g^{\perp}$, since $L\in
\mathbb{F}_{v}$, write $h=h_{1}\vee h_{2}\vee\cdots\vee h_{m}$,
where $h_{i}\wedge h_{j}=0$ for $i\neq j$ and each $h_{j}$ is
special. So each $h_{j}$ is a basic element. Since $h\wedge g=0$
implies that $h\wedge g_{i}=0$ for $i=1,2,\cdots,n$, so that
$h_{j}\wedge g_{i}=0$ for $j=1,2,\cdots,m$. From this we can
further obtain that $h_{j}\wedge a_{\lambda_{i}}=0$ for any
$\lambda_{i}\in \Lambda_{1}$ and $j=1,2,\cdots,m$. Using Claim 1,
each $h_{j}\in h^{\perp\perp}_{j}=a^{\perp\perp}_{\lambda_{k}}$
for some $\lambda_{k}\in\Lambda_{2}$ and $j=1,2,\cdots,m$.
So\vskip 4mm
\begin{center}
$h\in
\bigvee\limits_{j=1}^{m}h^{\perp\perp}_{j}=\bigvee\limits_{k=1}^{m}a^{\perp\perp}_{\lambda_{k}}
\subseteq
\bigvee\limits_{\lambda\in\Lambda_{2}}a^{\perp\perp}_{\lambda}$.\vskip
4mm
\end{center}
\noindent Thus
$g^{\perp}\subseteq\bigvee\limits_{\lambda\in\Lambda_{2}}a^{\perp\perp}_{\lambda}$.
Therefore
$g^{\perp}=\bigvee\limits_{\lambda\in\Lambda_{2}}a^{\perp\perp}_{\lambda}$.\vskip
3mm

Finally, since\vskip 2mm
\begin{center}
$L=\bigvee\limits_{\lambda\in\Lambda}a^{\perp\perp}_{\lambda}$,\vskip
2mm
\end{center}
\noindent using the results of Claim 2 and Claim 3, we then
have\vskip 2mm
\begin{center}
$L=\bigvee\limits_{\lambda\in\Lambda}a^{\perp\perp}_{\lambda}=
(\bigvee\limits_{\lambda_{i}\in\Lambda_{1}}a^{\perp\perp}_{\lambda_{i}})\vee
(\bigvee\limits_{\lambda\in\Lambda_{2}}a^{\perp\perp}_{\lambda})=g^{\perp\perp}\vee
g^{\perp}$.\vskip 2mm
\end{center}
\noindent Therefore $L\in \mathbb{T}$.$\hfill\Box$ \vskip 8mm

\noindent{\bf Acknowledgement }\vskip 4mm

This research was supported by the Science and Technology
Development Foundation of Nanjing University of Science and
Technology and the National Natural Science Foundation of China.
\vskip 8mm

\baselineskip 11pt \vskip 0.88 true cm

\end{document}